\documentclass[12pt, a4paper]{article}
\usepackage{mathtools,amsthm,amsmath,amsfonts,amssymb,tikz-cd,enumitem, graphicx,mathrsfs,bbm,stmaryrd,xcolor}

\tikzcdset{scale cd/.style={every label/.append style={scale=#1},
    cells={nodes={scale=#1}}}}
\usepackage{url}
\usepackage{esvect}
\usepackage{setspace}
\usepackage[utf8]{inputenc}
\usepackage{CJKutf8}
\counterwithin{equation}{section}
\usepackage[pdfencoding=auto,psdextra,colorlinks=true,linkcolor=blue,citecolor=blue,urlcolor = orange]{hyperref}
\usepackage[margin = 1.25in]{geometry}
\usepackage{kantlipsum}
\allowdisplaybreaks
\newtheorem{theorem}{Theorem}[section]
\newtheorem{definition}[theorem]{Definition}

\newtheorem{lemma}[theorem]{Lemma}

\newtheorem{remark}[theorem]{Remark}

\newtheorem{corollary}[theorem]{Corollary}

\newcommand{\citep}[1]{\cite{#1}}

\newcommand{\edmp}{\text{\normalfont End}}

\newcommand{\AuthorInfo}{}

\newcommand{\AddAuthor}[2]{%
 \appto{\AuthorInfo}{%
   \noindent\textbf{#1}\\#2\par\vspace{1em}
  }%
}

\AddAuthor{}{
Beijing International Center for Mathematical Research, Peking University\\
\texttt{\color{orange}hzhuang@pku.edu.cn}
}

\begin{document}
\title{Odd mapping cone Chern character}
\date{\today}
\author{Hao Zhuang}
\maketitle
\begin{abstract}
Our previous work established the mapping cone Chern character on the even-degree side, using the powers of the mapping cone curvature associated with the given mapping cone connection. In this paper, we carry out the odd-degree version using a trivial mapping cone connection and an automorphism on a trivial complex vector bundle. This version is based on the related odd-degree transgression. 
%

%
\end{abstract}

\tableofcontents
\section{Introduction}
Motivated by the $p$-filtered cohomology and the associated mapping cone model \cite{tanaka_tseng_2018, tty3rd, tty1st, tty2nd}, we apply the mapping cone connection \cite{tseng_and_zhou_symplectic_flat_connection_and_twisted_primitive2022, tseng_zhou_symplectic_flat_functional_characteristic_classes2022, tseng_and_zhou_2025mapping_yang_mills} to construct the associated transgression and characteristic classes \cite{transgression_primitive_2025, mappingconethomforms}. Those constructions are even-degree. As the classical transgression has an odd-degree version \cite{Cibotaru_Pereira_2022, Dai_Zhang_2006index_Toeplitz, getzler1993odd_Chern_character, wittendeformationweipingzhang, weipingzhangnewedition}, we hope to find an odd-degree transgression in the mapping cone cohomology.

In this paper, we build this odd-degree counterpart. We let $M$ be a smooth manifold, $\omega\in\Omega^2(M)$ be a closed smooth $2$-form, $E$ be a trivial rank $r$ complex vector bundle on $M$, $\delta$ be a trivial connection on $E$, and $z\in\mathbb{C}$ be a constant. 

\begin{definition}
\normalfont 
    We call 
\begin{align}
    A_z: \Omega^k(M,E)\oplus\Omega^{k-1}(M,E)&\to \Omega^{k+1}(M,E)\oplus\Omega^k(M,E) \nonumber\\
    (\alpha, \beta) &\mapsto (\delta\alpha + \omega\wedge\beta, z\alpha - \delta\beta)
\end{align}
a trivial mapping cone connection on $E$. Meanwhile, we call
\begin{align}
    \mathbf{tr}: \Omega^k(M,\edmp(E))\oplus\Omega^{k-1}(M,\edmp(E))&\to\Omega^k(M)\oplus\Omega^{k-1}(M) \nonumber\\
    (\gamma, \eta)&\mapsto (\text{tr}(\gamma), \text{tr}(\eta))
\end{align}
the trace map on pairs of $\edmp(E)$-valued forms. 
\end{definition}

Now, we let $\Gamma(\text{Aut}(E))$ be the space of smooth sections of $\text{Aut}(E)$. Given a family 
$g_t\in\Gamma(\text{Aut}(E))$ being smooth 
in terms of $t\in\mathbb{R}$, we have 
\begin{align}
    A_z(g_t, 0) = (\delta g_t, zg_t)\in\Omega^1(M,\edmp(E))\oplus\Omega^0(M,\edmp(E))
\end{align}
for each $t\in\mathbb{R}$.
Following \cite[Proposition 1.2]{getzler1993odd_Chern_character} and \cite[(1.50)]{wittendeformationweipingzhang} and using the product \cite[Section 2]{transgression_primitive_2025} of pairs of $\edmp(E)$-valued forms, we define: 
\begin{definition}\label{Chern character}
\normalfont
    For any odd number $2n+1$, we call 
    \begin{align}
        c_{2n+1}(g_t, A_z)\coloneqq \left(2\pi\sqrt{-1}\right)^{-n-1}\mathbf{tr}\left((g_t^{-1}, 0) A_z(g_t, 0)\right)^{2n+1}
    \end{align}
    the $(2n+1)$-th Chern form associated with $g_t$ and $A_z$. Then, we call 
    \begin{align}
        \text{ch}(g_t, A_z)\coloneqq \sum_{n\geqslant 0}\dfrac{n!}{(2n+1)!} c_{2n+1}(g_t, A_z) 
    \end{align}
    the Chern character form associated with $g_t$ and $A_z$. 
\end{definition}

Let $d$ be the de Rham differential. Then we have the mapping cone differential
\begin{align}
    d^\omega: \Omega^k(M)\oplus\Omega^{k-1}(M)&\to \Omega^{k+1}(M)\oplus\Omega^k(M) \nonumber\\
    (\alpha,\beta)&\mapsto (d\alpha + \omega\wedge\beta, -d\beta) 
\end{align}
satisfying $d^\omega\circ d^\omega = 0$. The following main result ensures that we indeed construct some characteristic classes in the mapping cone cohomology.
\begin{theorem}\label{main result}
For $n\in\mathbb{Z}_{\geqslant 0}$, we have 
    \begin{align}\label{closed form}
        d^\omega\mathbf{tr}\left((g_t^{-1}, 0) A_z(g_t, 0)\right)^{2n+1} = (0, 0)
    \end{align}
    and   
    \begin{align}\label{independent of t}
         & \dfrac{\partial}{\partial t}\mathbf{tr}\left((g_t^{-1}, 0) A_z(g_t, 0)\right)^{2n+1} \nonumber\\
        =\ & (2n+1) d^\omega\mathbf{tr}\left(\Big(g_t^{-1}\dfrac{\partial g_t}{\partial t}, 0\Big)\left((g_t^{-1}, 0) A_z(g_t, 0)\right)^{2n}\right). 
    \end{align}
\end{theorem}

By (\ref{closed form}) and (\ref{independent of t}), the mapping cone cohomology classes represented by Definition \ref{Chern character} are independent of $t$. The independence of these classes on $z$ is automatically revealed in the proof of Theorem \ref{main result}. In terms of other factors of dependence, we have the following corollary (cf. \cite[Corollary 1.18]{wittendeformationweipingzhang}). 
\begin{corollary}\label{corollary 1}
    Given any $g, h\in \Gamma(\mathrm{Aut}(E))$, we have
    \begin{align}
        & \mathbf{tr}\left(((gh)^{-1},0)A_z (gh, 0)\right)^{2n+1} \nonumber\\
        =\ & \mathbf{tr}\left((g^{-1},0)A_z (g, 0)\right)^{2n+1} + \mathbf{tr}\left((h^{-1},0)A_z (h, 0)\right)^{2n+1} + d^\omega\mu, 
    \end{align}
    where $\mu\in\Omega^{2n}(M)\oplus\Omega^{2n-1}(M)$. 
\end{corollary}

     In addition, like \cite[Corollary 1.19]{wittendeformationweipingzhang}, we have the  independence on the choice of the trivial connection $\delta$ on $E$. 

\begin{corollary}\label{corollary 2}
    Let $g\in \Gamma(\mathrm{Aut}(E))$. If $\widetilde{\delta}$ is another trivial connection on $E$, and $\widetilde{A}_z$ is the associated trivial mapping cone connection, then
    \begin{align}
        \mathbf{tr}\left((g^{-1},0)\widetilde{A}_z (g, 0)\right)^{2n+1}
    = \mathbf{tr}\left((g^{-1},0)A_z (g, 0)\right)^{2n+1} + d^\omega\nu 
    \end{align}
    with $\nu\in\Omega^{2n}(M)\oplus\Omega^{2n-1}(M)$. 
\end{corollary}

 In Section \ref{section 2}, we prove Theorem \ref{main result} and point out a key observation in Remark \ref{observation}. In Section \ref{section 3}, we prove Corollaries \ref{corollary 1}-\ref{corollary 2}. 

\section{Closedness and transgression}\label{section 2}
In this section, we prove Theorem \ref{main result}. We first clarify some objects. Suppose that 
\begin{align}
    e_1, e_2, \cdots, e_r
\end{align}
is a global frame of $E$, and that $\delta$ is associated with this frame. Then for any $f_1, \cdots, f_r\in C^\infty(M)$, we have 
\begin{align}
    \delta(f_1 e_1 + \cdots + f_r e_r) = df_1\otimes e_1 + \cdots + df_r\otimes e_r. 
\end{align}
Thus, if we write 
\begin{align}
    \begin{bmatrix}
        g_t(e_1), \cdots, g_t(e_r)
    \end{bmatrix} = \begin{bmatrix}
        e_1, \cdots, e_r
    \end{bmatrix}\begin{bmatrix}
        h_{11} & \cdots & h_{1r}\\
        \vdots &        & \vdots\\
        h_{r1} & \cdots & h_{rr}
    \end{bmatrix},
\end{align}
then $\delta g_t$ is given by 
\begin{align}\label{trivial connection local}
    \begin{bmatrix}
        \delta g_t(e_1), \cdots, \delta g_t(e_r)
    \end{bmatrix} = \begin{bmatrix}
        e_1, \cdots, e_r
    \end{bmatrix}\begin{bmatrix}
        dh_{11} & \cdots & dh_{1r}\\
        \vdots &        & \vdots\\
        dh_{r1} & \cdots & dh_{rr}
    \end{bmatrix}.
\end{align}

Now we recall the product of pairs of $\edmp(E)$-valued forms \cite[Section 2]{transgression_primitive_2025}: For $(\alpha_1, \beta_1) \in\Omega^k(M)\oplus\Omega^{k-1}(M)$, $(\alpha_2, \beta_2)\in\Omega^{\ell}(M)\oplus\Omega^{\ell-1}(M)$, and $G_1, H_1, G_2, H_2\in\Gamma(\edmp(E))$, we have the product
\begin{align}
    & (\alpha_1\otimes G_1, \beta_1\otimes H_1)(\alpha_2\otimes G_2, \beta_2\otimes H_2) \nonumber\\
    =\ & \left((\alpha_1\wedge\alpha_2)\otimes(G_1G_2), (\beta_1\wedge\alpha_2)\otimes(H_1G_2) + (-1)^k(\alpha_1\wedge\beta_2)\otimes (G_1H_2) \right).
\end{align} 
Thus, we see that 
\begin{align}
       &\left((g_t^{-1}, 0) A_z(g_t, 0)\right)^{2n+1} \nonumber\\
    =\ & \left(g_t^{-1}\delta g_t, g_t^{-1}zg_t\right)^{2n+1} \nonumber\\
    =\ & \left((g_t^{-1}\delta g_t)^{2n+1}, \sum_{k = 0}^{2n}(-1)^k (g_t^{-1}\delta g_t)^k (g_t^{-1}zg_t) (g_t^{-1}\delta g_t)^{2n-k}\right) \nonumber\\
    =\ & \left((g_t^{-1}\delta g_t)^{2n+1}, z\sum_{k = 0}^{2n}(-1)^k (g_t^{-1}\delta g_t)^{2n}\right)  \nonumber\\
    =\ & \left((g_t^{-1}\delta g_t)^{2n+1}, z(g_t^{-1}\delta g_t)^{2n}\right). 
\end{align}

Recall \cite[(1.40)]{wittendeformationweipingzhang}:
\begin{lemma}
For any even number $2n$, we have 
    \begin{align}\label{vanishing of trace in even power}
    \mathrm{tr}(g_t^{-1}\delta g_t)^{2n} = \dfrac{1}{2}\mathrm{tr}\left[(g_t^{-1}\delta g_t)^{2n-1}, g_t^{-1}\delta g_t\right] = 0. 
\end{align}
\end{lemma}

Thus, by (\ref{vanishing of trace in even power}), we immediately have 
\begin{align}\label{part 1}
    & \mathbf{tr}\left((g_t^{-1}, 0) A_z(g_t, 0)\right)^{2n+1} \nonumber\\
    =\ & \left(\text{tr}(g_t^{-1}\delta g_t)^{2n+1}, z\hspace{+0.4mm}\mathrm{tr}(g_t^{-1}\delta g_t)^{2n}\right) \nonumber\\
    =\ & \left(\text{tr}(g_t^{-1}\delta g_t)^{2n+1}, 0\right). 
\end{align}
In addition, we find  
    \begin{align}
        & \left(g_t^{-1}\dfrac{\partial g_t}{\partial t}, 0\right)\left((g_t^{-1}, 0) A_z(g_t, 0)\right)^{2n} \nonumber\\
       =\ & \left(g_t^{-1}\dfrac{\partial g_t}{\partial t}, 0\right)\left((g_t^{-1}\delta g_t)^{2n}, \sum_{k = 0}^{2n-1}(-1)^k (g_t^{-1}\delta g_t)^k (g_t^{-1}zg_t) (g_t^{-1}\delta g_t)^{2n-1-k}\right) \nonumber\\
       =\ & \left(g_t^{-1}\dfrac{\partial g_t}{\partial t}, 0\right)\left((g_t^{-1}\delta g_t)^{2n}, \sum_{k = 0}^{2n-1} z (-1)^k (g_t^{-1}\delta g_t)^{2n-1}\right) \nonumber\\
       =\ &  \left(g_t^{-1}\dfrac{\partial g_t}{\partial t}, 0\right)\left((g_t^{-1}\delta g_t)^{2n}, 0\right) \nonumber\\
       =\ & \left(g_t^{-1}\dfrac{\partial g_t}{\partial t}(g_t^{-1}\delta g_t)^{2n}, 0\right). 
    \end{align}
Thus, we have
\begin{align}\label{part 2}
    \mathbf{tr}\left(\Big(g_t^{-1}\dfrac{\partial g_t}{\partial t}, 0\Big)\left((g_t^{-1}, 0) A_z(g_t, 0)\right)^{2n}\right) = \left(\text{tr}\Big(g_t^{-1}\dfrac{\partial g_t}{\partial t}(g_t^{-1}\delta g_t)^{2n}\Big), 0 \right). 
\end{align}

Now we recall the classical closedness and transgression in the odd-degree case. The details can be found in \cite[Example 5.4]{Cibotaru_Pereira_2022}, \cite[Lemma 5.3]{Dai_Zhang_2006index_Toeplitz}, \cite[Propositions 1.2-1.3]{getzler1993odd_Chern_character}, \cite[(1.42), Lemma 1.17]{wittendeformationweipingzhang}, and \cite[(1.80), Lemma 1.7]{weipingzhangnewedition}.

\begin{theorem}[Getzler \cite{getzler1993odd_Chern_character}]
We have the $d$-closedness 
\begin{align}\label{closednessstatement}
    d\hspace{+0.5mm}\mathrm{tr}(g_t^{-1}\delta g_t)^{2n+1} = 0
\end{align}
and the odd-degree transgression 
    \begin{align}\label{part 3}
   \dfrac{\partial\mathrm{tr}(g_t^{-1}\delta g_t)^{2n+1}}{\partial t} = (2n+1)d\hspace{+0.5mm}\mathrm{tr}\Big(g_t^{-1}\dfrac{\partial g_t}{\partial t}(g_t^{-1}\delta g_t)^{2n}\Big). 
\end{align}
\end{theorem}
Thus, by (\ref{part 1}), (\ref{part 2}), (\ref{closednessstatement}), and (\ref{part 3}), we find that 
\begin{align}
    & d^\omega\mathbf{tr}\left((g_t^{-1}, 0) A_z(g_t, 0)\right)^{2n+1} \nonumber\\
    =\ & d^\omega\left(\text{tr}(g_t^{-1}\delta g_t)^{2n+1}, 0\right) \nonumber\\
    =\ & \left(d\hspace{+0.5mm}\mathrm{tr}(g_t^{-1}\delta g_t)^{2n+1}, 0\right) \nonumber\\
    =\ & (0,0), 
\end{align}
and that 
\begin{align}
    & \dfrac{\partial}{\partial t} \mathbf{tr}\left((g_t^{-1}, 0) A_z(g_t, 0)\right)^{2n+1} \nonumber\\ 
    =\ & \left(\dfrac{\partial\text{tr}(g_t^{-1}\delta g_t)^{2n+1}}{\partial t}, 0\right)   \nonumber\\
    =\ & \left((2n+1)d\hspace{+0.5mm}\mathrm{tr}\Big(g_t^{-1}\dfrac{\partial g_t}{\partial t}(g_t^{-1}\delta g_t)^{2n}\Big), 0\right)   \nonumber\\
    =\ & (2n+1)d^\omega\left(\mathrm{tr}\Big(g_t^{-1}\dfrac{\partial g_t}{\partial t}(g_t^{-1}\delta g_t)^{2n}\Big), 0\right)   \nonumber\\
    =\ & (2n+1) d^\omega \mathbf{tr}\left(\Big(g_t^{-1}\dfrac{\partial g_t}{\partial t}, 0\Big)\left((g_t^{-1}, 0) A_z(g_t, 0)\right)^{2n}\right). 
\end{align}
Theorem \ref{main result} is now proved. 

\begin{remark}\label{observation}
\normalfont
    Since the classical odd-degree transgression uses the trivial connection $\delta$, we may at first think that Theorem \ref{main result} can only be built for $z = 0$. However, the interesting thing is, due to (\ref{vanishing of trace in even power}), we can let $z$ be any complex number. This observation may be more important than the closedness and the transgression. The coincidence is, we do not need the even power of $g_t^{-1}\delta g_t$ to be $0$. Instead, we only need its trace to be $0$. In addition, considering the relation between the second entry in each pair of forms and the angular form on a circle bundle \cite[Appendix A]{tanaka_tseng_2018}, this observation may deserve a more detailed explanation. 
\end{remark}

\section{Factors of dependence}\label{section 3}
In this section, we prove Corollaries \ref{corollary 1} and \ref{corollary 2}. We recall the following fact \cite[Corollary 1.18]{wittendeformationweipingzhang}. 

\begin{lemma}
Given any $g, h\in \Gamma(\mathrm{Aut}(E))$, we have
    \begin{align}\label{homotopy}
        \mathrm{tr}((gh)^{-1}\delta(gh))^{2n+1} = \mathrm{tr}(g^{-1}\delta g)^{2n+1} + \mathrm{tr}(h^{-1}\delta h)^{2n+1} + d\gamma,
    \end{align}
    where $\gamma$ is a smooth $2n$-form on $M$. 
\end{lemma}

\begin{remark}\normalfont
    The proof of (\ref{homotopy}) has a K-theoretic interpretation \cite[Lemma 2.2.7]{atiyah2018k}. 
\end{remark}

We know the trivial connection $\delta$ is associated with a trivialization of $E$ (See (\ref{trivial connection local})). When we have another trivial connection, we recall \cite[Corollary 1.19]{wittendeformationweipingzhang}:   
\begin{lemma}
    Let $g\in \Gamma(\mathrm{Aut}(E))$. For another trivial connection $\widetilde{\delta}$ on $E$, we have 
    \begin{align}\label{change trivialization}
        \mathrm{tr}(g^{-1}\widetilde{\delta}g)^{2n+1} = \mathrm{tr}(g^{-1}\delta g)^{2n+1} + d\eta, 
    \end{align}
    where $\eta$ is a smooth $2n$-form on $M$. 
\end{lemma}

Therefore, by (\ref{homotopy}), we have 
\begin{align}
    & \mathbf{tr}\left(((gh)^{-1},0)A_z (gh, 0)\right)^{2n+1} \nonumber\\
    =\ & \left(\mathrm{tr}((gh)^{-1}\delta(gh))^{2n+1}, \mathrm{tr}((gh)^{-1}\delta(gh))^{2n}\right) \nonumber\\
    =\ & \left(\mathrm{tr}((gh)^{-1}\delta(gh))^{2n+1}, 0\right) \nonumber\\
    =\ & \left(\mathrm{tr}(g^{-1}\delta g)^{2n+1} + \mathrm{tr}(h^{-1}\delta h)^{2n+1} + d\gamma, 0\right) \nonumber\\
    =\ & \left(\mathrm{tr}(g^{-1}\delta g)^{2n+1}, 0\right) + \left(\mathrm{tr}(h^{-1}\delta h)^{2n+1}, 0\right) + d^\omega(\gamma, 0)  \nonumber\\
    =\ & \left(\mathrm{tr}(g^{-1}\delta g)^{2n+1}, \mathrm{tr}(g^{-1}\delta g)^{2n}\right) + \left(\mathrm{tr}(h^{-1}\delta h)^{2n+1}, \mathrm{tr}(h^{-1}\delta h)^{2n}\right) + d^\omega(\gamma, 0) \nonumber\\
    =\ & \mathbf{tr}\left((g^{-1},0)A_z(g,0)\right)^{2n+1} + \mathbf{tr}\left((h^{-1},0)A_z(h,0)\right)^{2n+1} + d^\omega(\gamma, 0).
\end{align}
Corollary \ref{corollary 1} is now proved. 

Finally, we let $\widetilde{A}_z$ be the trivial mapping cone connection associated with $\widetilde{\delta}$. By (\ref{change trivialization}), we have 
\begin{align}
    & \mathbf{tr}\left((g^{-1},0)\widetilde{A}_z (g, 0)\right)^{2n+1} \nonumber\\
    =\ & \left(\mathrm{tr}(g^{-1}\widetilde{\delta}g)^{2n+1}, 0\right) \nonumber\\
    =\ & \left(\mathrm{tr}(g^{-1}\delta g)^{2n+1}, 0\right) + (d\eta, 0) \nonumber\\
    =\ & \mathbf{tr}\left((g^{-1},0)A_z (g, 0)\right)^{2n+1} + d^\omega(\eta, 0).  
\end{align}
Corollary \ref{corollary 2} is now proved. 

\section*{Acknowledgments}
I thank Prof. Hongzhi Liu, Prof. Xiaobo Liu, Prof. Xiang Tang, Prof. Li-Sheng Tseng, Prof. Hang Wang, and Prof. Shanwen Wang for helpful discussions. Also, I thank Beijing International Center for Mathematical Research for the vibrating working environment. 

\bibliographystyle{abbrv}
\bibliography{mybib.bib}

@book{wittendeformationweipingzhang,
  title={Lectures on Chern-Weil Theory and Witten Deformations},
  author={Zhang, W.},
  isbn={9789812386588},
  year={2001},
  publisher={World Scientific},
  note = {\url{https://doi.org/10.1142/4756}}
}

@book{atiyah2018k,
  title={K-theory},
  author={Atiyah, M. F. },
  isbn={9780429962097},
  year={2018},
  publisher={CRC Press}, 
  note = {\url{https://doi.org/10.1201/9780429493546}}
}

@book{weipingzhangnewedition,
author = {Weiping Zhang and Huitao Feng},
 title = {Geometry and Analysis on Manifolds},
year = {2022},
edition = {{C}hinese},
publisher = {Higher Education Press},
note={\url{https://xuanshu.hep.com.cn/front/book/findBookDetails?bookId=60dedfb4adb85dae6a2f430b}}
}

@article{tanaka_tseng_2018,
    author = {Hiro Lee Tanaka and Li-Sheng Tseng},
    title = {Odd sphere bundles, symplectic manifolds, and their intersection theory},
    journal = {Cambridge Journal of Mathematics},
    volume = {6},
    number = {3},
    year = {2018},
    pages = {213--266},
    doi = {10.4310/CJM.2018.v6.n3.a1},
    issn = {2168-0930},
    publisher = {International Press of Boston},
    note = {\url{https://doi.org/10.4310/CJM.2018.v6.n3.a1}}
}

@article{tty1st,
author = {Li-Sheng Tseng and Shing-Tung Yau},
title = {{Cohomology and Hodge theory on symplectic manifolds: I}},
volume = {91},
journal = {Journal of Differential Geometry},
number = {3},
publisher = {Lehigh University},
pages = {383--416},
year = {2012},
doi = {10.4310/jdg/1349292670},
note = {\url{https://doi.org/10.4310/jdg/1349292670}}
}

@article{tty2nd,
author = {Li-Sheng Tseng and Shing-Tung Yau},
title = {{Cohomology and Hodge theory on symplectic manifolds: II}},
volume = {91},
journal = {Journal of Differential Geometry},
number = {3},
publisher = {Lehigh University},
pages = {417--443},
year = {2012},
doi = {10.4310/jdg/1349292671},
note = {\url{https://doi.org/10.4310/jdg/1349292671}}
}

@article{tty3rd,
author = {Chung-Jun Tsai and Li-Sheng Tseng and Shing-Tung Yau},
title = {{Cohomology and Hodge theory on symplectic manifolds: III}},
volume = {103},
journal = {Journal of Differential Geometry},
number = {1},
publisher = {Lehigh University},
pages = {83--143},
year = {2016},
doi = {10.4310/jdg/1460463564},
note = {\url{https://doi.org/10.4310/jdg/1460463564}}
}

@article{tseng_and_zhou_2025mapping_yang_mills,
  title={Mapping Cone Connections and their {Yang--Mills functional}},
  author={Tseng, Li-Sheng and Zhou, Jiawei},
  journal={Communications in Mathematical Physics},
  volume={406},
  number={7},
  pages={156},
  year={2025},
  publisher={Springer},
  note = {\url{https://doi.org/10.1007/s00220-025-05311-8}}
}

@article{tseng_and_zhou_symplectic_flat_connection_and_twisted_primitive2022,
  title={Symplectic flatness and twisted primitive cohomology},
  author={Tseng, Li-Sheng and Zhou, Jiawei},
  journal={The Journal of Geometric Analysis},
  volume={32},
  number={11},
  pages={282},
  year={2022},
  publisher={Springer},
  note = {\url{https://doi.org/10.1007/s12220-022-01018-7}}
}

@misc{tseng_zhou_symplectic_flat_functional_characteristic_classes2022,
  title={Symplectically flat connections and their functionals},
  author={Tseng, Li-Sheng and Zhou, Jiawei},
  year={2022},
  note = {\url{https://arxiv.org/abs/2210.03032}}
}

@misc{transgression_primitive_2025,
      title={Transgression in the primitive cohomology}, 
      author={Hao Zhuang},
      year={2025},
      eprint={2512.24920},
      archivePrefix={arXiv},
      primaryClass={math.DG},
      url={https://arxiv.org/abs/2512.24920}, 
      note = {\url{https://arxiv.org/abs/2512.24920}}
}

@article{mappingconethomforms,
      title={Mapping cone {T}hom forms}, 
      author={Hao Zhuang},
      year={2026},
      journal={Peking Mathematical Journal},
      volume={online published},
      note = {\url{https://doi.org/10.1007/s42543-026-00133-6}}
}

@article{Cibotaru_Pereira_2022, 
title={Nontame {Morse–Smale flows and odd Chern–Weil} theory}, 
volume={74}, 
DOI={10.4153/S0008414X21000353}, 
number={6}, 
journal={Canadian Journal of Mathematics}, 
author={Cibotaru, Daniel and Pereira, Wanderley}, 
year={2022}, 
pages={1579–1624}, 
note = {\url{https://doi.org/10.4153/S0008414X21000353}}
}

@article{getzler1993odd_Chern_character,
  title={The odd {Chern} character in cyclic homology and spectral flow},
  author={Getzler, Ezra},
  journal={Topology},
  volume={32},
  number={3},
  pages={489--507},
  year={1993}, 
  note = {\url{https://doi.org/10.1016/0040-9383(93)90002-D}}
}

@article{Dai_Zhang_2006index_Toeplitz,
  title={An index theorem for {Toeplitz} operators on odd-dimensional manifolds with boundary},
  author={Dai, Xianzhe and Zhang, Weiping},
  journal={Journal of Functional Analysis},
  volume={238},
  number={1},
  pages={1--26},
  year={2006},
  publisher={Elsevier}, 
  note = {\url{https://doi.org/10.1016/j.jfa.2006.05.005}}
}
\AuthorInfo
\end{document}